\documentclass{article}
\usepackage{amsmath,amstext,amssymb,amsopn,amsthm}

\theoremstyle{plain}
\newtheorem{thm}{Theorem}[section]
\newtheorem*{thm*}{Theorem}

\newtheorem{prop}[thm]{Proposition}
\newtheorem{cor}[thm]{Corollary}
\theoremstyle{definition}

\theoremstyle{remark}
\newtheorem*{rem*}{Remark}
\newtheorem{rem}{Remark}

\newcommand{\R}{\mathbb{R}}

\newcommand{\Z}{\mathbb{Z}}

\renewcommand{\geq}{\geqslant}

\newcommand{\pref}[1]{(\ref{#1})}

\DeclareMathOperator{\dist}{dist}

\def\a{\alpha}   
\def\o{\over}
\def\({\left(}
\def\){\right)}
\def\[{\left[}
\def\]{\right]}

\def\<{\langle}
\def\>{\rangle}
\def\A{{\bf A}^{(x,t)}}
\def\B{{\bf B}^{(x,t)}}
\def\as{$\alpha$-stable }

\title{Stable semigroups on homogeneous trees
  and hyperbolic spaces 
\footnotetext{2000 MS Classification:
    Primary 60J35;
    Secondary 47D03, 14M17. {\it Key words and phrases}: homogeneous trees,
    hyperbolic spaces, stable processes, heat kernels, metric spaces.
    Research 
    partially supported by MNI
    grant 1 P03A 020 28  and RTN Harmonic Analysis and Related Problems
    contract HPRN-CT-2001-00273-HARP}}
\author{ {\sc Andrzej St\'os}\\[5mm]  
Laboratoire de Math\'ematiques\\ Universit\'e Blaise Pascal\\
24, av. des Landais\\
63177 Aubi\`ere Cedex, France\\
\tt email: stos@math.univ-bpclermont.fr\\
and\\
Instytut Matematyki\\
 Politechnika Wroc\l{}awska\\
Wyb. Wyspia\'nskiego 27\\
50-370 Wroc\l{}aw, Poland
}
\date{}
\begin{document}

\maketitle

\vspace{1cm}

\abstract{We prove the kernel estimates related to subordinated
  semigroups on homogeneous trees. We study the long time 
  propagation problem. We exploit this to show exit time
  estimates for (large) balls. We use an abstract setting of metric
  measure spaces. This enables us to give these results for trees end
  hyperbolic spaces as well. Finally, we show some  estimates for the
  Poisson kernel of a ball.
}

\pagebreak

\section{Introduction}

In 1961 Getoor \cite{Ge} proposed subordinated semigroups in the
context of the real hyperbolic spaces.  It is only recently when the
corresponding kernel estimates were found (\cite{AJ}, \cite{GS}).

The aim of this paper is to give a corresponding result in the context
of homogeneous trees.  Our motivations come from the fact that such
structures make a discrete setting counterpart for hyperbolic spaces.
Large scale analogy holds not only in geometry but also in analysis,
see e.g. \cite{FTN}, \cite{FTP}, \cite{CMS}.

Our starting point is a diffusion semigroup considered in \cite{CMS}.
By subordination we obtain a new semigroup, which is referred to as to
the stable one. We show estimates for the corresponding 
kernel (Theorem \ref{ptx} below).  In the proof we use time-space
relations discovered in \cite{GS}.  On the other hand, our theorem
leads to a natural interpretation for the analogous result from
\cite{GS} (see remarks after the proof).

Next, we consider the long time propagation problem (Theorem
\ref{repart1}). It turns out that for large time $t$ the mass of our
kernel is distributed at distances comparable with $t^{2/\alpha}$.
We give two different proofs. First of them is of general nature and exploits 
some properties of the underlying diffusion semigroup.
This works
for hyperbolic spaces or Riemannian manifolds as well. The other
proof, a very simple one, shows that our Theorem \ref{ptx}
is useful as well.

Getoor \cite{Ge} raised the question of "stability" properties for
semigroups of this type.  Obviously, here we have neither classical
scaling, nor its {\em weak} form which is typical for e.g. fractals
\cite{BSS}.  However, one may interpret Theorem \ref{repart1} as an
{\em asymptotical} scaling property. A sample of its consequences is
given in the last section.

We conclude the paper by giving an application of Theorem \ref{ptx}.
We study exit time from balls for the {\it stable} process
corresponding to our semigroup. For related results we refer the
reader to \cite{GJ} or \cite{J}.  In general, we were inspired by the
approach from \cite{Ba1}, for stable case see \cite{BSS}.  The results
in section \ref{last} have their analogues in these papers. Observe,
however, that the argument of \cite{Ba1} and \cite{BSS} hinges on the
{\em Ahlfors-regularity} of the measure, i.e. polynomial volume
growth.  Clearly, this excludes the homogeneous trees and hyperbolic
spaces.  Our contribution is to make it available for stable processes
in exponential volume growth setting.  Moreover, we give a proof in an
abstract framework of metric measures spaces (cf. \cite{G}). The
interplay between \pref{vol} and \pref{ptx2} below may be of
independent interest.  In this way, we get our results for homogeneous
trees and hyperbolic spaces at the same time.

Finally, using the Ikeda-Watanabe formula we give estimates for the
Poisson kernel for balls.

\section{Preliminaries}
Consider the nearest-neighbor Laplacian $\Delta$ and the related heat
semigroup ${\cal H}_t$ with continuous time on a homogeneous tree $X$
of degree $q+1$ with $q\ge 2$, i.e.
\begin{equation*}
  \Delta f(x) = f(x) -\frac{1}{q+1}\sum_{y\sim x}f(y), \quad x\in X
  \quad {\rm and} \quad {\cal  H}_t = e^{-t\Delta}, \quad t>0.
\end{equation*}
See \cite{CMS} for detailed exposition. We adopt the general setting
from this paper. For the reader's convenience we recall definitions
needed in what follows.  In particular, let $h_t$ denote the
corresponding heat kernel and $h_t^\Z$ the heat kernel in the
one-dimensional case.  Moreover, set $\gamma={2\sqrt{q} \o q+1}$ so
that $b_2=1-\gamma$ is the bottom of the spectrum of the Laplacian
acting on $L^2(X)$.

We adopt the convention that $c$ (without subscripts) denotes a
generic constant whose value may change from one place to another. To
avoid some curiosities occasionally we write $\tilde c$, $c'$...  with
the same properties. Numbered constants (with subscripts) always keep
their particular value throughout the current theorem or proof.  We
often write $f\asymp g$ to indicate that there exists $c>0$ such that
$c^{-1}<f/g<c$. Similarly, $f(x)\asymp g(x)$, $x\to\infty$, means $f
\asymp g$ for $x$ large enough.

The kernel $h_t$ is known to satisfy the following estimates
\cite{CMS}:
\begin{equation*}
  h_t(x) \asymp {e^{-b_2 t} \o t}\phi_0(x) h_{t\gamma}^\Z(|x|+1),
\end{equation*}
where
\begin{equation}\label{spher}
  \phi_0(x) = \left( 1+ {q-1\o q+1} |x| \right) q^{-{|x|\o2}}, \quad x\in X
\end{equation}
is the {\em spherical function},
\begin{equation*}
  h_t^\Z(j) = e^{-t}I_{|j|}(t),\quad t>0,\, j\in\Z,
\end{equation*}
and $I_\nu(t)$ stands for the modified Bessel function of the first
kind.  Consequently,
\begin{equation}  \label{hk}
  h_t(x) \asymp {e^{-t} \o t} \phi_0(x) I_{1+|x|}(t\gamma), \quad
  t>0,\,\,x\in X.
\end{equation}
In what follows we fix $\alpha\in(0,2)$ and consider the {\it
  subordinate semigroup }
\begin{equation*}
  e^{-t \Delta^{\alpha/2}} = \int_0^\infty e^{-u\Delta}\eta_t(u)du,
\end{equation*}
where the {\em subordinator} $\eta_t(\cdot)$ is a function (defined on
$\R^+$) determined by its Laplace transform,
\begin{equation*}
  {\cal L}[\eta_t(\cdot)](\lambda) = e^{-t\lambda^{\alpha/2}}.
\end{equation*}
For the corresponding kernels we have
\begin{equation}\label{sub}
  p_t(x) = \int_0^\infty h_u(x)\eta_t(u) du.
\end{equation}
Sometimes we refer to $p_t(x)$ as to the \as kernel.  For more
details concerning this construction we refer the reader e.g.  to
\cite{Be}.

\section{\as  kernel}
The main result may be stated as follows.
\begin{thm}\label{ptx}
  For any constants $K,M > 0$
  \begin{equation}\label{ptx1}
    p_t(x) \asymp \left\{
      \begin{array}{ll}
        \phi_0(x)t^{-3/ 2}\exp(-t(1-\gamma)^{\alpha/2}), 
        & |x|<K t^{1/2},\, t\ge 1,
        \\
        \phi_0(x)t|x|^{-2-{\alpha/2}}q^{-{|x|/2}}, & |x|>Mt^{2/\alpha}>0.
      \end{array}\right.
  \end{equation}
\end{thm}

\begin{proof}
  First, we collect some auxiliary estimates for Bessel function
  $I_\nu(z)$. Recall its integral representation (e.g. \cite{GR},
  (8.431.1))
  \begin{equation*}
    I_\nu(z)= {(z/2)^{\nu}\o \Gamma(\nu+{1/2})\sqrt{\pi}}
    \int_{-1}^1 (1-u^2)^{\nu-1/2}e^{-zu}du
    ={(2\pi z)^{-1/2}e^z \o 2^{\nu-1/2}\Gamma(\nu+1/2)}
    \int_0^{2z} [ u(2-u/z) ]^{\nu-1/2}e^{-u}du.
  \end{equation*}
  Clearly, the last integral is bounded above by
  $2^{\nu-1/2}\Gamma(\nu+1/2)$ so that
  \begin{equation}\label{ileq}
    I_\nu(z) \le cz^{-1/2}e^z, \quad z>0,\, \nu>0.
  \end{equation}
  Let us recall that (\cite{CMS})
  \begin{equation}
    \label{iequiv}
    I_\nu(z) \asymp {e^{\sqrt{\nu^2+z^2}}\o \sqrt{z+\nu}} 
    \left( {z \o \nu+\sqrt{\nu^2+z^2}  }\right)^\nu, 
    \quad \nu\ge 1,\, z>0.
  \end{equation}
  Assume that $z>\max(1,a\nu^2)$ with some $a\in (0,1)$ and $\nu>1$.
  Thus,$\sqrt{\nu^2+z^2}-z\le a/2$ so that $\exp(\sqrt{\nu^2+z^2})
  \asymp \exp(z)$ (in the lower bound there is a constant that depends
  on $a$).  Clearly, $\sqrt{z+v}\asymp\sqrt{z}$ and the quotient in
  the parentheses in \pref{iequiv} is bounded above by 1. Moreover,
  \begin{equation*}
    {z^\nu \o (\nu+\sqrt{\nu^2+z^2})^\nu} 
    \ge
    { 1\o ( \sqrt{a} / \sqrt{z}+\sqrt{1+a/z} )^{\sqrt{az}} }
    \ge
    { 1\o \left( 1+{2\sqrt{a}/ \sqrt{z} }
      \right)^{\sqrt{z}/(2\sqrt{a})\times 2a} } \ge {1\o e^{2a}}.
  \end{equation*}
  
  Consequently, we obtain the desired simplification
  \begin{equation}\label{ieq}
    I_\nu(z) \asymp z^{-1/2}e^{z}, \quad z>\max(1,a\nu^2),\, \nu\ge 1.
  \end{equation}
  
  We recall the exact estimates of the densities $ \eta_t(\cdot)$
  which will be fundamental in what follows (see e.g. \cite{GS}). We
  have
  \begin{equation}
    \label{eta1}
    \eta_t(u)\asymp t^\frac{1}{2-\a} u^{- \frac{4-\a}{4-2\a}} \exp \left(
      -c_1t^\frac{2}{2-\a} u^{-\frac{ \a}{ 2-\a}} \right), \qquad
    t^{-2/\alpha}u<c,
  \end{equation}
  where $\displaystyle c_1=c_1(\alpha)=\frac{2-\a}{2}\left(
      \frac{\a}{2}
    \right)^\frac{ \a}{ 2-\a}$ and
  \begin{equation}
    \label{eta2}
    \eta_t(u)\asymp tu^{-1-\alpha/2}, \qquad t^{-2/\alpha}u>c.
  \end{equation}
  
  According to \pref{eta1} and \pref{eta2}, it is convenient to split
  the integral \pref{hk} as follows
  \begin{eqnarray}
    \nonumber
    p_t(x) & =&  \int_0^{c_0t^{2/\alpha}}h_u(x)\eta_t(u)du +
    \int_{c_0t^{2/\alpha}}^\infty h_u(x)\eta_t(u)du 
    \\ 
    \nonumber & =& 
    \phi_0(x)t^{1\o 2-\alpha} \int_0^{c_0t^{2/\alpha}}
    e^{-u}I_{1+|x|}(\gamma u)u^{-{4-\alpha\o 4-2\alpha}-1} 
    \exp(-c_1t^{2\o2-\alpha}u^{-{\alpha\o 2-\alpha}}) du
    \\ \label{decomp} & + &
    \phi_0(x)t \int_{c_0t^{2/\alpha}}^\infty e^{-u}I_{1+|x|}(\gamma u)
    u^{-2-{\alpha/2}}du
    \\ \nonumber &=&
    \phi_0(x)\( \A+\B\).
  \end{eqnarray}
  
  Now, we assume that $c_0=1$ and $|x|\le K\sqrt{t}$ with $x$ and $t$
  large enough. Note that neither $x$, nor $t$ is fixed.  It follows
  that $(1+|x|)^2\le(1+K\sqrt{t})^2 \le \gamma t^{2/\alpha}$.  Hence,
  by \pref{ieq} with $a=1$ we get
  \begin{equation*}
    I_{1+|x|}(\gamma u) \le c u^{-1/2}e^{\gamma u},\quad u>t^{2/\alpha}.
  \end{equation*}
  In consequence,
  \begin{eqnarray*}
    \B &\le & 
    ct\int_{t^{2/\alpha}}^\infty e^{-(1-\gamma)u}u^{-(5+\alpha)/2} du
    \\& \le & 
    c t^{-{5/\alpha}}
    \int_{t^{2/\alpha}}^\infty e^{-(1-\gamma)u}du 
    \\ & = &
    ct^{-5/\alpha}
    e^{-(1-\gamma)t^{2/\alpha}}
  \end{eqnarray*}
  
  To estimate $\A$ let us split it as follows
  \begin{equation*}
    \A = t^{1\o 2-\alpha} \( \int_0^{ {\alpha t/2} } +
    \int_{{\alpha t /2}}^{t^{2/\alpha}} \)
    e^{-u}I_{1+|x|}(\gamma u)u^{-{4-\alpha\o 4-2\alpha}-1} 
    \exp(-c_1t^{2\o2-\alpha}u^{-{\alpha\o 2-\alpha}}) du
    = \A_1+\A_2.
  \end{equation*}
  Now, apply \pref{ieq} to the integral $\A_2$. After simple change of
  the variable $u\to tu$, we get
  \begin{eqnarray*}
    \A_2 &\asymp& c t^{-1}\int_{\alpha\o2}^{t^{{2\o\alpha}-1}}
    u^{-{4-\alpha\o4-2\alpha}-{3\o2}}\exp(-t ((1-\gamma)u+c_1
    u^{-{\alpha\o2-\alpha}})) du.
  \end{eqnarray*}
  Observe that the minimum of function $ p(u)= (1-\gamma) u
  +c_1u^{-{\alpha\o2-\alpha}}$ is attained at
  \begin{equation*}
    u_0={ \(\alpha c_1\o 2-\alpha\)^{1-{\alpha\o2}} \o
      (1-\gamma)^{1-{\alpha /2}}}. 
  \end{equation*}
  Since
  \begin{equation*}
    \(\alpha c_1 \o 2-\alpha \) ^{1-{\alpha\o 2}}
    =  \left[
      {\alpha\o 2-\alpha}{2-\alpha\o 2}
      \left(\alpha\o 2 \right)^{\alpha\o 2-\alpha}
    \right]^{2-\alpha\o 2} 
    \nonumber
    = \left[
      \left(\alpha\o 2\right)^{2-\alpha\o2}
      \left(\alpha\o 2\right)^{\alpha\o 2}
    \right] 
    \nonumber
    = {\alpha\o 2},
  \end{equation*}
  we get
  \begin{equation*}
    u_0  = {\alpha/2 \o (1-\gamma)^{1-{\alpha/2}}}.
  \end{equation*}
  Hence, for $t$ large enough $u_0$ is in the integration range and
  $p(u_0) = (1-\gamma)^{\alpha/2}$.  Obviously, our integral is
  bounded by integrals with limits fixed
  \begin{equation*}
    \int_{\alpha\o2}^{u_0}  \le \int_{\alpha\o2}^{t^{ {2\o\alpha}-1}}
    \le \int_0^\infty.
  \end{equation*}
  The Laplace method \cite{O} applied to the extreme members of this
  inequality gives the same result, so that we obtain the asymptotic
  of our integral:
  \begin{equation*}
    ct^{-1/2}e^{-tp(u_0)}, \quad  t\to \infty.
  \end{equation*}
  Consequently,
  \begin{equation*}
    \A_2 \asymp t^{-{3/2}}\exp(-(1-\gamma)^{\alpha/2}t), \quad
    |x|<K\sqrt{t}
  \end{equation*}
  and $t\ge1$, say.  Similarly, using \pref{ileq} we get
  \begin{equation*}
    \A_1 \le  c t^{-1}\int_{\alpha\o2}^{t^{{2\o\alpha}-1}}
    u^{-{4-\alpha\o4-2\alpha}-{3\o2}}\exp(-t ((1-\gamma)u+c_1
    u^{-{\alpha\o2-\alpha}})) du.
  \end{equation*}
  Since the minimum of $p(u)$ is {\em not} attained in $(0,\alpha/2)$,
  in this case the Laplace method gives the following lower bound:
  \begin{equation*}
    \A_1\le  c t^{-2}\exp(-p(\alpha/2)t).
  \end{equation*}
  It follows that $p_t(x)\asymp \A_2$ and the first of the desired
  estimates follows.
  
  Now, assume that $|x|>Mt^{2/\alpha}$. Since we consider large $|x|$
  only (or even $|x| \to \infty$), we may and do put $|x|-1$ in place
  of $|x|$ when estimating $p_t(\cdot)$. This simplifies the notation.
  We put $c_0=aM$ in the decomposition \pref{decomp}, where
  $a\in(0,1)$ is to be specified later. Then, by \pref{iequiv} and the
  elementary inequalities $e^{\sqrt{|x|^2+\gamma^2u^2}}\le
  e^{|x|+\gamma u}$, $|x|+\sqrt{|x|^2+\gamma^2u^2}\ge 2|x|$, we get
  \begin{eqnarray*}
    \A & = & 
    t^{1\o 2-\alpha}\int_0^{aMt^{2/\alpha}}
    e^{-u}I_{|x|}(\gamma u)u^{-{4-\alpha\o 4-2\alpha}-1} 
    e^{-c_1t^{2\o2-\alpha}u^{-{\alpha\o 2-\alpha}}} du
    \\   &\le & 
    c|x|^{\alpha\o 4-2\alpha} \int_0^{a|x|}
    { e^{\sqrt{|x|^2+\gamma^2u^2}-u} (\gamma u)^{|x|}
      u^{-{4-\alpha\o4-2\alpha}-1} 
      e^{-c_1t^{2\o2-\alpha}u^{-{\alpha\o 2-\alpha}}} \o
      \sqrt{|x|+\gamma u} \left(|x|+\sqrt{|x|^2+\gamma^2 u^2}
      \right)^{|x|} } du 
    \\ & \le &
    c|x|^{{\alpha\o 4-2\alpha}-{1\o2}}\(ae\gamma\o2\)^{|x|}
    \int_0^{a|x|} 
    e^{-(1-\gamma)u} u^{-{4-\alpha\o 4-2\alpha}-1}
    e^{-c_1t^{2\o2-\alpha}u^{-{\alpha\o 2-\alpha}}}du.
  \end{eqnarray*}
  Clearly, the last integral is convergent and bounded above by a
  constant independent of $|x|$. Therefore,
  \begin{equation}\label{atx}
    \A\le  c|x|^{{\alpha\o 4-2\alpha}-{1\o2}}\(ae\gamma\o2\)^{|x|}.
  \end{equation}
  On the other hand, again by \pref{iequiv} and the change of variable
  $u\to ux$, we obtain
  \begin{eqnarray*}
    \B & = & 
    t \int_{aM t^{2/\alpha}}^\infty I_{|x|}(\gamma u)
    u^{-2-{\alpha/2}}e^{-u}du 
    \\ & \ge &
    c t \int_{a|x|}^\infty 
    {e^{\sqrt{|x|^2+\gamma^2u^2}-u} \o \sqrt{|x|+\gamma u} }
    { (\gamma u)^{|x|} u^{-2-{\alpha/2}} \o
      \left(|x|+\sqrt{|x|^2+\gamma^2u^2}\right)^{|x|} }du
    \\ & \ge & 
    c{t \gamma^{|x|} |x|^{-{\alpha+3\o2}} }
    \int_{a}^\infty {e^{|x|(\sqrt{1+\gamma^2u^2}-u)} u^{-2-{\alpha/2}+|x|}
      \o \sqrt{1+\gamma u} \left(1+\sqrt{1+\gamma^2u^2}\right)^{|x|}} du
    \\ & \asymp &
    {t\gamma^{|x|} |x|^{-{\alpha+3\o2}}} \int_{a}^\infty
    e^{|x|\(\sqrt{1+\gamma^2u^2}-u+\log(u)-\log(1+\sqrt{1+\gamma^2u^2})\)}
    {u^{-2-{\alpha /2}} \o \sqrt{1+\gamma u}} du.
  \end{eqnarray*}
  Observe that the same calculation with the lower limit of
  integration equal to 0 gives the opposite bound. Let
  \begin{equation*}
    p(u)=\sqrt{1+\gamma^2u^2}-u+\log(u)-\log(1+\sqrt{1+\gamma^2u^2})
  \end{equation*}
  and $g=\sqrt{1+\gamma^2u^2}$. Then $p'(u)=
  -1+g/u$ and, consequently, $p(u)$ attains the maximum at $u_0={q+1\o
    q-1}>1$. Hence, $u_0$ belongs to the integration range for
  integrals in both upper and lower bound for $\B$. Consequently, by
  the Laplace method, both of them have the same asymptotic as $|x|\to
  \infty$. Since
  \begin{equation*} p(u_0)=\sqrt{1+{4q\o(q-1)^2}}-{q+1\o
      q-1}+\log\(q+1\o (q-1)\(1+\sqrt{1+{4q\o(q-1)^2}}\) \) = -\log\(
    2q\o q+1 \) = -\log(\gamma\sqrt{q}),
  \end{equation*}
  it follows that
  \begin{equation*}
    \B \asymp t |x|^{-2-{\alpha/2}} e^{|x|(\log \gamma -\log(\gamma\sqrt{q}))}
    = t |x|^{-2-{\alpha/2}}q^{-{|x|/2}}, \quad|x|\ge
    Mt^{2/\alpha}
  \end{equation*}
  and $|x|$ is large enough (and hence for $|x|>1$).  Moreover, if we
  take $a=1/e$ then ${ae\gamma / 2} \le q^{-1/2} $ so that $\A=o(\B)$,
  $|x|\to \infty$ and $p_t(x) \asymp \B$.  The assertion follows.
\end{proof}

\begin{rem}
  Our theorem can be compared with the following result of \cite{GS}.
  For reader's convenience we give it below, specialized to the (real)
  hyperbolic space $\mathbb{H}^n$.  The corresponding \as 
  kernel and spherical function are denoted with the tilde.

  \begin{thm*}\label{gs}[\cite{GS}, Corollary 5.6]
    Let $|\rho|=(n-1)/2$.  For any constants $K,M>0$ and $t+|x|>1$ we
    have
    \begin{equation}\label{hs}
      \tilde p_t(x) \asymp \left\{
        \begin{array}{ll}
          \tilde \phi_0(x)\,t^{-{3/2}}\,e^{-|\rho|^\a t}, & 
          |x|\le K\,t^{1/2} 
          \\ \tilde \phi_0(x)\, t|x| ^{-2-{\alpha/2}}
          e^{-|\rho||x|}, & |x| \ge M\,t^{2/\alpha}.
        \end{array} \right. 
    \end{equation}
  \end{thm*}
  
  In the context of hyperbolic space (or, more generally, symmetric
  space of non-compact type), the parameter $|\rho|$ plays a double
  role: it is the square root of the bottom of the spectrum of the
  Laplace-Beltrami operator; at the same time, the volume growth of
  the ball of the radius $r$ is equivalent to $e^{2|\rho| r}$,
  $r\to\infty$.  One may ask, whether it is the spectral data or the
  geometry which appears in the above estimates. The comparison with
  Theorem \ref{ptx} gives us a natural interpretation: in the first
  part (i.e. in the long time asymptotics) we deal with the spectral
  data, in the other case the volume growth intervenes.
\end{rem}
\begin{rem} 
  Note that  for the
  remaining region $ Kt^{1/2}< |x| < M t^{2/\alpha}$, in the
  continuous setting there
  is no {\em simple explicit } estimate of $\tilde p_t(x)$ (see
  \cite{GS}, Corollary 5.6). 
\end{rem}

The Brownian motion and \as processes in $\R^d$ share the
same type of long time heat repartition. Namely, with the standard
understanding that $\alpha=2$ corresponds to the Brownian motion, for
$A_1<A_2$ we have
\begin{equation*}
  \int_{A_1t^{1/\alpha}\le |x|\le A_2t^{1/\alpha}}
  p_t(x)dx = c(A_1,A_2) \in(0,1).
\end{equation*}
This follows immediately from the scaling property
\begin{equation} \label{scaling}
  p_t(x)=t^{-d/\alpha}p_1(t^{-1/\alpha}x).
\end{equation}
Moreover, $c(A_1,A_2)\to 1$ if $A_1\to 0$ and $A_2\to \infty$ so that
\begin{equation}\label{excl}
  \int_{ A_1t^{\beta}\le |x|\le A_2t^{\beta}}
  p_t(x)dx \to 0,\quad t\to \infty,
\end{equation}
provided $\beta\not=1/\alpha$ (cf. \cite{AS}, p. 50).

On the other hand, for the Brownian motion in the (real) hyperbolic
space $\mathbb{H}^n$, a non-classical phenomenon of concentration was
observed in \cite{D}.  Namely,
\begin{equation*}
  \int_{A_1 t \le    |x| \le A_2 t}h_t(x)dx \to 1,\quad t\to \infty,
\end{equation*}
provided $A_1 < n-1 <A_2$.  The change of the space-time scaling
should be noted. This result was sharpened and generalized to
symmetric space setting (\cite{AS}, \cite{Ba}). In the context of
homogeneous trees the analogous result was shown in \cite{MS} and
\cite{W}:
\begin{equation*}
  \sum_{R_0t-r(t)\le |x| \le R_0t+r(t)} h_t(x)\to 1,\quad t\to\infty,
\end{equation*}
where $R_0=(q-1)/(q+1)$ and $r(t)$ is a positive function such that
$r(t)t^{-1/2}\to \infty$, $t\to\infty$.  This might suggest a
hypothesis of the same kind for our kernel $p_t(x)$, e.g.
\begin{equation*}
  \sum_{A_1 t^{2/\alpha}\le |x| \le A_2t^{2/\alpha}} p_t(x)
  \to 1,\quad t\to\infty.
\end{equation*}
The following theorem shows that this is not the case.

\begin{thm}\label{repart1}
  For $0<A_1<A_2$ let $R(t)=\{ (x,t)\in X\times \R^+: A_1 t^{2/\alpha}
  \le |x| \le A_2t^{2/\alpha} \}$. Then there exist $c_1$ and $c_2$
  such that
  \begin{equation}\label{asymp4}
    0<c_1<\sum_{x\in R(t) } p_t(x) <c_2<1, 
    \quad t\to\infty.
  \end{equation}
  Conversely, for any given $0<c_1<1$ ($0<c_2<1$ resp.) there exist
  $A_1$ and $A_2$ such that \pref{asymp4} holds true with some $c_2$
  ($c_1$ resp.).
\end{thm}
\begin{proof}
  Set $R_0=(q-1)/(q+1)$ and let $R_1$, $R_2$ be such that
  $R_1<R_0<R_2$. Then, by Theorem 1 of \cite{MS}, we have
  \begin{equation}\label{conv}
    \sum_{R_1u \le |x|\le R_2u}h_u(x) \to 1,\qquad u\to \infty.
  \end{equation}
  Moreover, let $c_3=A_1/R_1$ and $c_4=A_2/R_2$. We require
  additionally that $R_1$ and $R_2$ be close to $R_0$ so that
  $c_3<c_4$.  Then $c_3t^{2/\alpha}< u < c_4t^{2/\alpha}$ yields
  \begin{equation}\label{conv1}
    |x|\in (R_1u,R_2u) \implies x\in R(t).
  \end{equation}
  From the definition of $p_t(x)$, \pref{conv} and \pref{conv1}, we
  get
  \begin{eqnarray*}
    \sum_{x\in R(t)} p_t(x) 
    & = & 
    \int_0^\infty  \(\sum_{x\in R(t)} h_u(x) \) \eta_t(u)du
    \\ &\ge & 
    \int_{c_3t^{2/\alpha}}^{c_4t^{2/\alpha}}  
    \( \sum_{R_1u \le |x| \le R_2u} h_u(x) \) \eta_t(u) du
    \\ &\to  & 
    \int_{c_3t^{2/\alpha}}^{c_4t^{2/\alpha}} \eta_t(u) du, \quad t\to
    \infty. 
  \end{eqnarray*}
  Formally, the last integral depends on $t$. By the scaling property
  \pref{scaling}, however, it evaluates to
  \begin{equation}\label{eta}
    t^{-2/\alpha}
    \int_{c_3t^{2/\alpha}}^{c_4t^{2/\alpha}}\eta_1(t^{-2/\alpha}u)du = 
    \int_{c_3}^{c_4}\eta_1(u)du = c_0.
  \end{equation}
  This is an absolute constant which depends on $c_3$, $c_4$ and
  $\alpha$ only. The lower bound in the first assertion follows.
  Since the lower bound is true for {\em any} $A_1<A_2$, the mass of
  the annulus $R(t)$ (with $A_1$ and $A_2$ fixed) is strictly less
  than 1.  In other words, $c_2<1$ in \pref{asymp4} and there is no
  {\em total} mass concentration observed.  The proof of the first
  assertion is complete.
  
  The second assertion follows from the fact that $c_0$ in \pref{eta}
  can be required to take any value in $ (0,1)$. Indeed, if $A_1\to 0$
  and $A_2\to \infty$ then we may fix $R_1<R_0<R_2$ independently of
  $A_1$ and $A_2$, so that $c_3\to 0$ and $c_4\to \infty$. Since
  $\int_0^\infty \eta_t(u) = 1$ we may require $c_1$ to be arbitrarily
  close to 1.
  
  Further, fix any $0<\tilde c_2$. The upper bound for the mass of the
  annulus $\tilde R(t)=\{ (x,t)\in X\times \R^+: \tilde A_1
  t^{2/\alpha} \le |x| \le \tilde A_2t^{2/\alpha} \}$ follows from the
  lower bound for $R(t)$ provided $A_2< \tilde A_1$. Since the mass of
  $R(t)$ can be required to be arbitrarily close to 1, $t\to\infty$,
  the mass of $\tilde R(t)$ is (asymptotically) smaller than $\tilde
  c_2$.  The proof is complete.
\end{proof}
The following corollary is an analogue of the classical counterpart
\pref{excl}.
\begin{cor}
  For $0<\tilde A_1< \tilde A_2$ and some $\beta>0$ let $\bar R(t)=\{
  (x,t)\in X\times \R^+: \tilde A_1 t^{\beta} \le |x| \le \tilde
  A_2t^{\beta} \}$.  If $\beta\not= 2/\alpha$ then
  \begin{equation*}
    \sum_{x\in \bar R(t) } p_t(x) \longrightarrow 0,\quad t\to \infty.
  \end{equation*}
\end{cor}
\begin{proof}
  For $t$ large enough, $R(t)$ and $\bar R(t)$ are disjoint.
\end{proof}

Evidently, space-time scaling in \pref{asymp4} is characteristic for
the Brownian motion in hyperbolic spaces and homogeneous trees. On the
other hand, the concentration phenomenon is not observed. From the
probabilistic point of view this may be explained by the influence of
the long jumps of the corresponding stable process.  Indeed, the {\em
  L\'evy measure} is of the same exponential order as volume growth
because it arises from the second estimate in \pref{ptx1}. Actually,
we have
\begin{cor}
  Let $\nu(x):= \lim_{t\to 0} p_t(x)/t$ be the L\'evy measure for our
  semigroup. Then
  \begin{equation*}
    \nu(x) \asymp |x|^{-1-\alpha/2}q^{-|x|}, \quad |x|\ge 1.
  \end{equation*}
\end{cor}
\begin{proof}
  From Theorem \ref{ptx} and \pref{spher} we get
  \begin{equation*}
    \nu(x) \asymp \phi_0(x)|x|^{-2-\alpha/2}q^{-|x|/2} \asymp 
    |x|^{-1-\alpha/2}q^{-|x|}.
  \end{equation*}
\end{proof}

Clearly, the proof of Theorem \ref{repart1} with minor modifications
only can be applied in the context of the symmetric spaces with
Theorem 1 of \cite{AS} instead of \pref{conv}.  We prefer, however, to
take the opportunity given by Theorem 2 of that article to state our
result in the more general setting of manifolds.  For reader's
convenience, we recall the framework.  We assume that $M$ is a
complete, noncompact Riemannian manifold with the volume growth
controlled by
\begin{equation*}
  \mathrm{vol}(B(x,r)) = O(r^\kappa e^{2Kr}), \quad r\to\infty,
\end{equation*}
with some constants $\kappa$ and $K$, and the spectral gap $E^2= \inf
\mathrm{spec}(-\Delta)>0 $.  In general we have $E\le K$, while in
symmetric spaces $E=K=|\rho|$.  Set $R_1=2(K-\sqrt{K^2-E^2})$,
$R_2=2(K+\sqrt{K^2-E^2})$.  Let $A(t)$ be a function such that
\begin{equation*}
  \begin{array}{lcl}
    \displaystyle A(t) - \frac{\kappa-1}{2\sqrt{K^2-E^2}}\log t \nearrow \infty
    &  \mathrm{if} & K<E,\\
    A(t) = (2\kappa t\log t)^{1/2}
    & \mathrm{if} & K=E \mathrm{ ~~and~~ } \kappa>0,\\
    A(t)t^{-1/2}\nearrow \infty
    & \mathrm{if} & K=E \mathrm{ ~~and~~ } \kappa=0.
  \end{array}
\end{equation*}
Since the heat kernel depends on two variables (and is denoted by
$h_t(x,y)$), we fix arbitrary $y\in M$ and redefine slighlty $R(t) =
\{ (x,t)\in X\times \R^+: A_1 t^{2/\alpha} \le d(x,y) \le
A_2t^{2/\alpha} \}$.  By Theorem 2 from \cite{AS}
\begin{equation*}
  \int_{R_1t-A(t)\le d(x,y)\le R_2t+A(t) }h_t(x,y)\to 1,\quad t\to \infty.
\end{equation*}
Note that in any case we may and do require $A(t)=o(t)$, which is
essential for our proof to work (cf. \pref{conv1}).  Thus, we arrive
at
\begin{cor}\label{repart2}
  For any constants $0<A_1<A_2$ there exist $c_1$ and $c_2$ such that
  \begin{equation*}
    0<c_1 < \int_{R(t)}p_t(x)dx < c_2 <1,\quad t\to \infty.
  \end{equation*}
  Changing $A_1$ and $A_2$ we may require $c_1$ to be close to 1 or
  $c_2$ to be close to $0$.
\end{cor}

Below we include an alternative approach that relies directly on the
\as kernel estimates \pref{ptx1}.  It shows, in a sense, that
the \as  kernel mass covered by Theorem \ref{ptx} is large
enough to be useful in some applications.

\begin{proof}[Second proof of Theorem \ref{repart1}]
  For $x\in R(t)$ we have
  \begin{equation}\label{estim}
    p_t(x) \asymp t\phi_0(x)|x|^{-2-{\alpha/2}}q^{-{|x|/2}}.
  \end{equation}
  By \pref{spher},
  \begin{equation*}
    \phi_0(x) \asymp |x|q^{-{|x|/2}}, \quad |x|\to \infty.
  \end{equation*}
  Therefore,
  \begin{equation*}
    \sum_{x\in R(t)} p_t(x)
    \asymp t \sum_{x\in R(t) } |x|^{-1-\alpha/ 2} q^{-|x|} .
  \end{equation*}
  Since the function under the consideration depends only on the
  distance, we use ``polar coordinates''.  At each sphere $\{ |x|=n
  \}$ we have exactly $(q+1)q^n$ vertices, so that
  \begin{equation*}
    \sum_{x\in R(t)} p_t(x) 
    \asymp 
    t \sum_{A_1 t^{2/\alpha}\le n \le A_2 t^{2/\alpha}}
    n^{-1-{\alpha/2}}.
  \end{equation*}
  Since obviously
  \begin{equation*}
    \int_{a-2}^{b+2}y^{-1-\alpha/2}dy \asymp \int_a^b y^{-1-\alpha/2}dy
  \end{equation*}
  when $a\to \infty$ and $a/b<c_0<1$, it follows that
  \begin{equation*}
    \sum_{x\in R(t)} p_t(x)    
    \asymp t \int_{A_1 t^{2/\alpha}}^{A_2t^{2/\alpha}} y^{-1-{\alpha/2}}dy.
  \end{equation*}
  Clearly, the last integral behaves as
  \begin{equation}\label{conv2}  
    (A_1^{-\alpha/2}-A_2^{-\alpha/2}) (\alpha/2)^{-1}t^{-1}, \quad 
    t\to\infty.  
  \end{equation}
  The lower bound in \pref{asymp4} follows.  To obtain inequality
  $c_2<1$ in the upper bound, it is enough to take $A_1$ sufficiently
  large. To allow any value of $A_1<A_2$, we can repeat here the
  argument following \pref{eta} in the previous proof. The assertion
  follows.
\end{proof}
\begin{rem*}
  Clearly, \pref{estim} holds also for $\bar R(t)$ with
  $\beta>2/\alpha$ as well. In this case, \pref{conv2} implies that
  \begin{equation*}
    \sum_{x\in \bar R(t)}p_t(x)\to 0,\quad t\to \infty.
  \end{equation*}
  However, this direct argument fails for $\beta<2/\alpha$.  Actually,
  if \pref{estim} holded for $|x|\ge At^\beta$ with some
  $\beta<2/\alpha$, then we would obtain $t^{-\alpha\beta/2}$ in
  \pref{conv2}. Consequently, the mass of the annulus goes to
  infinity, which is impossible (cf. also Proposition 5.4 of
  \cite{GS}, where this was done for the line $|x|=t^\beta$ by an
  independent argument).
\end{rem*}

\section{Exit time}\label{last}
We conclude our work by giving an application of Theorem \ref{ptx}.
Since the results below are very similar for both homogeneous trees
and hyperbolic spaces, we are tempted to use the following notation of
metric spaces.

Let $(E,\mu)$ be a metric space with a measure $\mu$ that supports a
heat kernel $p_t(x,y)$ in the sense of the axiomatic definition 2.1 of
\cite{G}.  For the reader's convenience we recall it shortly. We
assume that $p_t(\cdot,\cdot)$ is a $\mu\times \mu$ nonnegative
measurable function and for $\mu$-almost all $x,y\in E$ and all $s,t
>0$ we have $p_t(x,y) = p_t(y,x)$,
\begin{equation*}
  \int_E p_t(x,y)d\mu(y)  =1, \quad p_{t+s}(x,y) =\int_E
  p_t(x,z)p_s(z,y)d\mu(z),
\end{equation*}
and for each $u\in L^2(E,\mu)$
\begin{equation*}
  \int_E p_t(x,y)u(y)d\mu(y) \stackrel{L^2}{\longrightarrow} u(x),
  \quad t\to 0^+. 
\end{equation*}

In the case of the hyperbolic spaces or homogeneous trees we have
$p_t(x,y)= p_t(d(x,y))$, where $d(x,y)$ is the distance.  Under some
general additional assumptions on $X$, this kernel gives rise to the
associated Markov process $(X_t, P_x)$, i.e.
\begin{equation*}
  P_x[X_t\in B] = \int_B p_t(x,y)d\mu(y).
\end{equation*}
For simplicity, we suppose that the space is homogeneous, i.e. there
exists a function $V(r)$, called a volume growth, such that $V(r) =
\mu(B(x,r))$, $x\in E$.  It can be seen that for the proofs below this
assumption is not essential and we could deal with non-homogeneous
version $V(x,r)$ as well.

Moreover, assume that there exist $A\ge 1$ and $c_1<1$ such that
\begin{equation}\label{vol}
  V(r) \le c_1V(r+A) \quad \text{and} \quad V(r+1)\asymp V(r), \quad r\ge 1.
\end{equation}
Actually, this covers the case of trees and hyperbolic spaces (with
e.g. $A=1$).

Furthermore, assume that for any $M>0$
\begin{equation}\label{ptx2}
  p_t(x,y) \asymp td(x,y)^{-1-\alpha/2}V(d(x,y))^{-1}, \quad
  d(x,y)>Mt^{2/\alpha}, \quad d(x,y)>1. 
\end{equation} 
This is clearly satisfied in the context of trees and hyperbolic
spaces as well (cf. Theorem \ref{ptx} and \pref{hs} resp.).

Note that the first part of \pref{vol} implies that
$\lim_{r\to\infty}V(r)=\infty$.  In particular, our space is not
contained in any ball. Below we use this fact without further mention.

\begin{prop} \label{exit1} For any $M>0$ and $r>1$ we have
  \begin{equation*}
    P_x[X_t\notin B(x,r)] \asymp tr^{-\alpha/2}, \quad r>Mt^{2/\alpha}.
  \end{equation*}
\end{prop}
\begin{proof}
  By \pref{ptx2} we get
  \begin{eqnarray*}
    P_x[X_t\notin B(x,r)] &\asymp& \int_{d(x,y)>r}p_t(x,y)d\mu(y)\\
    & \asymp  & t\sum_{k=0}^\infty \int_{r+k<d(x,y)\le r+k+1}
    d(x,y)^{-1-\alpha/2} V(d(x,y))^{-1}d\mu(y) \\
    &\le &
    c t\sum_{k=0}^\infty
    (r+k)^{-1-\alpha/2} V(r+k)^{-1}(V(r+k+1)-V(r+k)).
  \end{eqnarray*}
  Clearly, by \pref{vol} we get
  \begin{equation*}
    V(r+k)^{-1}(V(r+k+1)-V(r+k))  \le c.
  \end{equation*}
  Moreover, by a comparison of the series with the corresponding
  integral it can be easily seen that
  \begin{equation*}
    \sum_{k=0}^\infty    (r+k)^{-1-\alpha/2} 
    = r^{-1-\alpha/2} + \sum_{k=1}^\infty (r+k)^{-1-\alpha/2} 
    \le r^{-\alpha/2} + \int_r^\infty z^{-1-\alpha/2}dz 
    \le cr^{-\alpha/2}
  \end{equation*}
  and the upper bound in the assertion follows.
  
  On the other hand we have similarly
  \begin{eqnarray*}
    P_x[X_t\notin B(x,r)] &\asymp  & 
    t\sum_{k=0}^\infty \int_{r+kA<d(x,y)\le r+(k+1)A}
    d(x,y)^{-1-\alpha/2} V(d(x,y))^{-1}d\mu(y) \\
    &\ge &
    c t\sum_{k=0}^\infty
    (r+kA+A)^{-1-\alpha/2} \frac{V(r+kA+A)-V(r+kA)}{V(r+kA+A)}.
  \end{eqnarray*}
  Again, by \pref{vol}
  \begin{equation*}
    \frac{V(r+kA+A)-V(r+kA)}{V(r+kA+A)}
    = 1-\frac{V(r+kA)}{V(r+kA+A)} \ge 1-c_1>0.
  \end{equation*}
  Moreover,
  \begin{equation*}
    \sum_{k=0}^\infty (r+kA+A)^{-1-\alpha/2}
    \ge \int_{r+A}^\infty z^{-1-\alpha/2}dz
    = c(r+A)^{-\alpha/2} \ge cr^{-\alpha/2},
  \end{equation*}
  since $r>1$.  The proof is complete.
\end{proof}
For a measurable set $D$ define the exit time $\tau_D = \inf\{t\ge
0;\; X_t\notin D\}$. Then
\begin{prop}\label{exit2}
  For any $M>0$ and $r>1$ we have
  \begin{equation*}
    P_x[\tau_{B(x,r)} < t ] \le ctr^{-\alpha/2}, \quad r>Mt^{2/\alpha}.
  \end{equation*}
\end{prop}
\begin{proof}
  The proof follows the lines of \cite{Ba1} (or \cite{BSS}).  Since it
  is short, we sketch it for the reader's convenience.  Denote
  $T=\tau_{B(x,2r)}$. Then
  \begin{eqnarray*}
    P_x[T<t] & = &
    P_x[X_t\notin B(x,r);\, T<t] + P_x[X_t\in B(x,r);\, T<t]\\
    & \le &  P_x[X_t\notin B(x,r)] + P_x[X_t\in B(x,r);\, T<t] = A+B.
  \end{eqnarray*}
  By Proposition \ref{exit1} we obtain $A\le c tr^{-\alpha/2}$.  By
  the strong Markov property we have
  \begin{eqnarray*}
    B & = & E_x[P_{X(T)}[ X_{t-u}\in B(x,r)]_{|u=T};\, T<t]\\
    & \le &
    \sup_{u\le t}\sup_{z\in B(x,2r)^c} E_x[P_z[X_u\in B(x,r)];\, T<t]\\
    & \le & \sup_{u \le t}\sup_{z\in B(x, r)^c}
    E_x[P_z[X_u \notin B(z,r) ];\, T<t]\\
    & \le &  ctr^{-\alpha /2}.
  \end{eqnarray*}
  The proof is complete.
\end{proof}
\begin{thm} \label{exit3}
  For $r>1$
  \begin{equation*}
    E_y\tau_{B(x,r)} \le  cr^{\alpha/2}, \quad y\in B(x,r)
  \end{equation*}
  and
  \begin{equation*}
    E_x\tau_{B(x,r)} \asymp  r^{\alpha/2}.
  \end{equation*}
\end{thm}
\begin{proof}
  For any $y\in B(x,r)$ by Proposition \ref{exit1} we have
  \begin{equation*}
    P_y[\tau_{B(x,r)} >t]  \le  
    P_y[X_t\in B(x,r) ] 
    =  1-P_y[X_t\notin B(x,r)] 
    \le  1-ctr^{-\alpha/2}
  \end{equation*}
  provided that $r>Mt^{2/\alpha}$ with some $M>0$. Let
  $t_0=r^{\alpha/2}$ so that for some $c_0$ we get
  \begin{equation}\label{tau}
    P_y[\tau_{B(x,r)} > t_0]  \le  1-c_0.
  \end{equation}
  Then, by Markov property, for $k=1,2,...$ we have
  \begin{eqnarray*}
    P_y[ \tau_{B(x,r)} > (k+1) t_0 ] & = &
    P_y [ \tau_{B(x,r)}\circ\theta_{t_0} > k t_0,\; \tau_{B(x,r)}>t_0 ] \\
    & = & E_y[ P_{X(t_0)} [ \tau_{B(x,r)}>kt_0 ] ;\; \tau_{B(x,r)}>t_0 ] \\
    & \le & P_y[\tau_{B(x,r)}>t_0]\sup_{z \in B(x,r)   }
    P_z[\tau_{B(x,r)}>kt_0]
  \end{eqnarray*}
  (here $\theta$ stands for the standard shift operator on the space
  of trajectories).  By induction we get
  \begin{equation*}
    P_y[\tau_{B(x,r)}>kt_0] \le (1-c_0)^k, \qquad y\in B(x,r),\;\;k=0,1,2,...
  \end{equation*}
  Thus,
  \begin{equation*}
    E_y\tau_{B(x,r)}  =  \int_0^\infty P_y[\tau_{B(x,r)}>t]dt 
    \le 
    \sum_{k=0}^\infty t_0P_y[\tau_{B(x,r)}>kt_0]
    \le r^{\alpha/2}\sum_{k=0}^\infty (1-c_0)^k
  \end{equation*}
  and the upper bound in the asssertion follows.
  
  On the other hand, let $t_1=c_1r^{\alpha/2}$ with $c_1$ to be
  specified below.  From Proposition \ref{exit2} we get
  \begin{equation*}
    P_x[\tau_{B(x,r)} < t_1] \le c_1 c_2.
  \end{equation*}
  Observe that the constant $c_2$ above does not depend on $c_1$
  provided $c_1<1$. Hence, we may and do choose $c_1$ small enough to
  get $c_1c_2<1$.  It follows that
  \begin{equation*}
    E_x\tau_{B(x,r)} \geq t_1 P_x[\tau_{B(x,r)}>t_1] \geq
    (1-c_1c_2)t_1 \asymp  r^{\alpha/2}
  \end{equation*}
  The proof is complete.
\end{proof}

\section{Poisson kernel}
In this section we give estimates for the Poisson kernel for balls.
Since in general it  follows ideas of \cite{BSS}, we give only a short
sketch of the construction.  For more 
detailed exposition we refer the reader to  sections 5
and 6 of that article.  
Since the results in what follows are similar for both the
homogeneous trees and hyperbolic spaces, we continue to use the
notation introduced in the previous section. It should be noted,
however, that this concernes the results only. The details of proofs
in this section should be verified separately for each geometrical
setting.

In what follows, we assume that for $x,y\in X$ the following limit
exists
\begin{equation*}
  N(x,y) =  \lim_{t\to 0} \frac{p_t(x,y)}{t}>0.
\end{equation*}
This is verified whenever our \as kernel arises by a
subordination of a reasonable diffusion with $\eta_t$ described above.
Clearly, the case of homogeneous trees and hyperbolic spaces is
included. From \pref{ptx2} it follows that
\begin{equation}\label{nxy}
  N(x,y) \asymp d(x,y)^{-1-\alpha/2}V(d(x,y))^{-1}, \quad d(x,y)\ge 1. 
\end{equation}
Let
\begin{equation}\label{nxy1}
  n(x,E) = \int_E N(x,y)d\mu(y).
\end{equation}
For an open set $D$ let $(P_t^D)$ be the semigroup generated by the
process killed on exiting $D$, i.e.
\begin{equation*}
  P_t^Df(x) = E_x[f(X_t);\, t<\tau_D].
\end{equation*}
This semigroup possesses transition densities denoted by $p^D_t(x,y)$
(see \cite{ChZ}; the argument applies here as well). Let $G_D(x,y)$ be
the Green function for $D$, i.e. the potential for $(P_t^D)$:
\begin{equation*}
  G_D(x,y) = \int_0^\infty p^D_t(x,y)dt.
\end{equation*}

With these definitions one verifies the assumptions of the
Ikeda-Watanabe formula (see \cite{IW} or \cite{BSS}). For homogeneous
trees and hyperbolic spaces this is straightforward and we omit the
details. We note, however, that at this point each geometrical
structure is examined separately. We get
\begin{prop}[Ikeda-Watanabe formula]
  Assume that $D\subset X$ is an open nonempty bounded set, $E\subset
  X$ is a Borel set and $\dist(D,E)>0$.Then
  \begin{equation*}
    P_x[X_{\tau_D}\in E]  = \int_D G_D(x,y) n(x,E) d\mu(y).
  \end{equation*}
\end{prop}
In particular, by \pref{nxy1} we get that $P_x[X_{\tau_D}\in \cdot]$
is absolutely continuous w.r. to $\mu$ on $(\bar D)^c$ (this is
meaningful for the hyperbolic spaces only). Let $P_D(x,\cdot)$ denote
its density (i.e. Poisson kernel).
\begin{prop} 
  For any $x_0\in X$ and $r\ge 1$ let $D=B(x_0,r)$.  Then
  \begin{equation*}
    P_D(x,z) \le c
    \frac{r^{\alpha/2}V(2r)}{d(x,z)^{1+\alpha/2}V(d(x,z))}, 
    \quad z\in B(x_0,3r)^c,\,\,x\in D.
  \end{equation*}
  If $r\ge 2$ then
  \begin{equation*}
    P_D(x,z) \ge  c
    \frac{r^{\alpha/2}}{V(2r)d(x,z)^{1+\alpha/2}V(d(x,z))},
    \quad z\in D^c,\,\, x\in B(x_0,r/2).
  \end{equation*}
\end{prop}
\begin{proof} By \pref{nxy} we have 
  \begin{equation}\label{pdxy}
    P_D(x,z) \asymp  \int_D
    \frac{G_D(x,y)}{d(y,z)^{1+\alpha/2}V(d(y,z))} d\mu(y).
  \end{equation}
  Clearly, $d(y,z) \asymp d(x,z)$. Moreover, for the hyperbolic spaces
  and homogeneous trees we have $V(r) \asymp C_1^r$ where $C_1$
  depends on the dimension or the degree, respectively. It follows
  that
  \begin{equation*}
    V(d(y,z)) \ge V(d(x,z) - d(x,y)) \ge V(d(x,z)-2r) \asymp
    V(2r)^{-1}V(d(x,y)). 
  \end{equation*}
  Since $\displaystyle \int_D G_D(x,y)d\mu(y) = E_x\tau_D$ the upper
  bound in the assertion follows by Theorem \ref{exit3}.
  
  On the other hand, fix $x\in B(x_0,r/2)$. Then $d(y,z) \le c
  d(x,z)$, $y\in D$, $z\in D^c$. Similarly as before $V(d(y,z)) \le
  V(d(y,x)+d(x,z)) \asymp V(2r)V(d(x,z))$.  Moreover, $E_x\tau_D \ge
  E_x \tau_{B(x,r/2)} \asymp r^{\alpha/2}$.  By \pref{pdxy} the lower
  bound follows. We are done.
\end{proof}

\section*{Acknowledgements}
The hospitality of the MAPMO laboratory during author's postdoc stay
is gratefully acknowledged.


\begin{thebibliography}{99}
  \bibliographystyle{plain}
  
\bibitem{AJ} J.-P. Anker, L. Ji, {\it Heat kernel and Green function
    estimates on noncompact symmetric spaces}, Geom. Funct.  Anal
  9(6)(1999), 1035-1091.
  
\bibitem{AS} J.-P. Anker, A. G. Setti, {\it Asymptotic finite
    propagation speed for heat diffusion on certain Riemannian
    manifolds}, J. Funct. Analysis, 103(1)(1992), 50-61.
  
\bibitem{BG} R.M. Blumenthal, R.K. Getoor, {\it Some theorems on
    stable processes}, Trans. Amer. Math. Soc., 95(1960), 263-273.
  
\bibitem{Ba} M. Babillot, {\it A probabilistic approach to heat
    diffusion on symmetric spaces}, J. Theor. Probability 7(1994),
  599-607.
  
\bibitem{Ba1} M. T. Barlow, {\it Diffusion on fractals}, in: {\it
    Lectures on Probability Theory and Statistics, Ecole d'Ete de
    Probabilites de Saint-Flour XXV - 1995}, Lecture Notes in
  Mathematics no. 1690, Springer-Verlag, New York 1999, 1-121.

  
\bibitem{Be} J. Bertoin, {\it L\'evy Processes}, Cambridge University
  Press, Cambridge 1996.
  
\bibitem{BSS} K. Bogdan, A. St\'os, P. Sztonyk, {\it Harnack
    inequality for stable processes on $d$-sets}, Studia Math.
  158(2)(2003), 163-198.
  
\bibitem{CMS} M. Cowling, S. Meda, A. G. Setti, {\it Estimates for
    functions of the Laplace operator on homogeneous trees}, Trans.
  Amer. Math. Soc. 352(2000), 4271-4293.
  
\bibitem{ChZ} K. L. Chung, Z. Zhao, {\it From Brownian motion to
    Schr\"odinger's equation}, Springer-Verlag, New York 1995.
  
\bibitem{D} E.B. Davies, {\it Heat kernels and spectral theory},
  Cambridge University Press, 1989.
  
\bibitem{E} Erdelyi et al., eds., {\it Tables of integral transforms},
  II, McGraw-Hill, New York, 1954.
  
\bibitem{FTN} A. Fig\`a Talamanca, C. Nebia, {\em Harmonic analysis
    and representation theory for goups acting on homogeneous trees},
  London Math. Soc. Lecture Notes Series, n. 162, Cambridge University
  Press, 1991.
  
\bibitem{FTP} A. Fig\`a Talamanca, M. Picardello, {\em Harmonic
    analysis on free groups}, Lecture Notes in Pure and Applied
  Mathematics, n. 87, Marcel Dekker, 1983.
  
\bibitem{G} A. Grigoryan, {\it Heat kernel and function theory on
    metric measure spaces}, dans {\em Heat kernels and Analysis on
    manifolds, graphs and metric spaces, Paris 2002} Contemporary
  Math. 338(2003), 143-172.  cit\'e {\bf [1]}
  
\bibitem{GR} I. S. Gradstein, I. M. Ryzhik, {\it Table of integrals,
    series and products}, $6^{th}$ edition, Academic Press (London),
  2000.
  
\bibitem{Ge} R. K. Getoor, {\it Infinitely divisible probabilities on
    the hyperbolic plane}, Pacific J. Math.  11(1961), 1287-1308.
  
\bibitem{GJ} P. Graczyk, T. Jakubowski, {\em Exit times and Poisson
    kernels of the Ornstein-Uhlenbeck diffusion}, preprint.
  
\bibitem{GS} P. Graczyk, A. St\'os, {\it Transition density estimates
    for stable processes on symmetric spaces}, Pacific J.  Math.
  217(2004), no. 1, 87-100.
  
\bibitem{IW} N. Ikeda, S. Watanabe, {\it On some relations between the
    harmonic measure and the Lévy measure for a certain class of
    Markov processes}, J. Math. Kyoto Univ., 2(1962), 79-95.
  
\bibitem{J} T.Jakubowski, {\em Estimates of exit times for stable
    Ornstein-Uhlenbeck processes}, preprint.
  
\bibitem{MS} G. Medolla, A.G. Setti, {\it Long time heat diffusion on
    homogeneous trees}, Proc. Amer. Math. Soc., 128(6)(1999),
  1733-1742.
  
\bibitem{O} F. Olver, {\it Asymptotics and special functions},
  Academic Press, New York - London, 1974.
  
\bibitem{W} Woess, {\it Heat diffusion on homogeneous trees (note on a
    paper by Medolla and Setti)}, Bollettino U. M. I. 8 4-B (2001),
  703 - 709.

\end{thebibliography}
\end{document}